\documentclass[12pt]{article}
\title{Complex determinantal processes and $H^1$ noise}
\date{\today}
\author{Brian Rider
\ and B\'alint Vir\'ag\
}

\usepackage{amsfonts}
\usepackage{graphicx}
\usepackage{amsmath}
\usepackage{amsthm}
\usepackage{amssymb}

\theoremstyle{theorem}
    \newtheorem{theorem}{Theorem}
    \newtheorem{lemma}[theorem]{Lemma}
    \newtheorem{proposition}[theorem]{Proposition}
    \newtheorem{corollary}[theorem]{Corollary}

\theoremstyle{definition} 

\theoremstyle{remark} 

\newcommand\be{\begin{equation}}

\newcommand\ee{\end{equation}}

\newcommand{\comment}[1]{}

\newcommand{\UU}{{\mathbb U}}
\newcommand{\R}{{\mathbb R}}
\newcommand{\reals}{{\mathbb R}}
\newcommand{\CC}{{\mathbb C}}
\renewcommand{\SS}{{\mathbb S}}

\newcommand{\ev}{\mbox{\bf E}}

\newcommand{\one}{{\mathbf 1}}

\newcommand{\dist}{\rm dist}
\newcommand{\disti}{\mbox{\rm dist}_\iota}
\newcommand{\supp}{\mbox{\rm supp}}

\newcommand{\Var}{\mbox{\rm Var}}

\newcommand{\Cov}{\mbox{\rm Cov}}
\newcommand{\Cum}{\mbox{\rm Cum}}
\newcommand{\sm}{{\raise0.3ex\hbox{$\scriptstyle \setminus$}}}

\newcommand{\re}[1]{(\ref{#1})}

\newcommand{\ep}{\varepsilon}
\newcommand{\ra}{\rightarrow}

\newcommand{\KK}{\tilde K}

\newcommand{\eps}{\varepsilon}

\newcommand{\Ki}{K^\iota}

\newcommand{\inv}{\nu}
\newcommand{\nablai}{\nabla_\iota}

\begin{document}
\maketitle
\begin{abstract}
For the plane, sphere, and hyperbolic plane we consider the
canonical invariant determinantal point processes $\mathcal
Z_\rho$ with intensity $\rho d\nu$, where $\inv$ is the
corresponding invariant measure. We show that as $\rho\to\infty$,
after centering, these processes converge to invariant $H^1$
noise. More precisely, for all functions $f\in H^1(\inv) \cap
L^1(\inv)$ the distribution of $\sum_{z\in \mathcal Z_\rho}
f(z)-\frac{\rho}{\pi} \int f d \inv$ converges to Gaussian with
mean 0 and variance $ \frac{1}{4 \pi} \|f\|_{H^1}^2$.
\end{abstract}


\begin{figure}[h]
\centering
\includegraphics[height=1.4in]{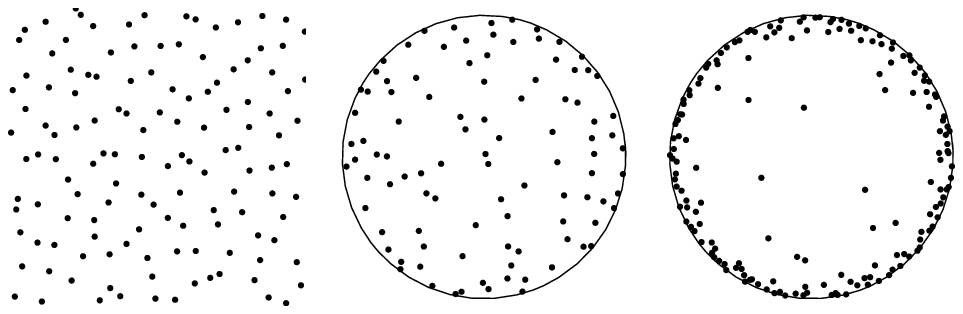}
\end{figure}

\section{Introduction}

Determinantal processes are point processes with a built-in pairwise
repulsion. They were first considered by Macchi \cite{Mac} as a
model for fermions in quantum mechanics, and have since been
understood to arise naturally in a number of contexts, from
eigenvalues of random matrices to random spanning trees and
non-intersecting paths, see
 \cite{Borodin,Burt,Diac,HKPV,Lyons,LyonsSt,PerVir,Shirai,Sosh}.

A point process $\mathcal Z$ on $\CC$
is {\bf determinantal}  if for disjoint sets $D_1, \dots, D_k$ we
have \be \label{joint}
 \ev \Bigl[    \prod_{i = 1}^k  \#(D_i \cap \mathcal Z)  \Bigr]  =
 \int_{\Lambda^k}  \det \Bigl(  K(x_i, x_j) \Bigr)_{1 \le i,j \le k}
 d \mu(x_1) \cdots d \mu (x_k),
 \ee for each $k \ge 1$. Here $K(x,y)$ is a hermitian symmetric
measurable function and $\mu$ is some {\bf reference measure}. The
integrand is often called the {\bf joint intensity} or correlation
function of the point process.

Conversely, if such a function $K$ defines a self-adjoint integral
operator $\cal K$ on $L^2( \Lambda, \mu)$ which is locally trace
class with all eigenvalues lying in $[0, 1]$, then there exists a
point process satisfying \re{joint}. In this case we speak of the
determinantal process $(K, \mu)$. In a weak sense, the points
repel one another because the determinant vanishes on the
diagonal.

The processes we consider are defined by two properties. First, they
correspond to a projection $K$ to a subspace of analytic functions
with respect to a radially symmetric reference measure. Second,
their distribution is invariant under the symmetries of their
underlying space $\Lambda$. The latter is either the complex plane
$\CC$, 2-sphere $\SS$, or hyperbolic plane $\UU$. We will think of
them as a subsets of $\CC$, or strictly speaking $\CC\cup\{\infty \}$,
though for the sphere $\SS$ usually there is no harm in ignoring the
point at infinity.

These properties are uniquely satisfied by a family of processes
indexed by a single density parameter $\rho>0$ on each space, see
Krishnapur \cite[Theorem 3.0.5]{Manju}, who discovered several
remarkable properties of these processes.

\vspace{.2cm}

\noindent {\bf Planar model.}  Here $\Lambda = \CC$ and for any
$\rho > 0$,  consider the kernel
\[
   {\check K}_{\rho}(z,  w) =   e^{\rho z \bar w}  \mbox{ with respect to }
   d \mu_{\rho} (z)  = \frac{\rho}{\pi} e^{-\rho |z|^2} d z.
\]
Here, as in the sequel,  $dz$ stands for Lebesgue area measure on
$\CC$.

Note that the kernel ${\check K}_{\rho}$ is the projection onto
the span of the orthonormal set $\{ \sqrt{\frac{\rho^k}{k!}} z^k
\}_{k=0}^{\infty}$ in $L^2(\CC , \mu_\rho)$, and so the above pair
defines a determinantal process of infinitely many points in the
complex plane, see, for example \cite{HKPV}.

 \vspace{.2cm}

\noindent {\bf Spherical model.}  The space $\Lambda$ is $\SS=\CC
\cup \infty$, the two-sphere. Now $\rho$ is integer valued, $\rho =
1,2,\dots$ and the $({\check K}, \mu)$ pair reads
\[
   {\check K}_{\rho}(z,  w) =  \sum_{k =0}^{\rho -1} { \rho -1 \choose k} (z {\bar w})^{k}=(1+z{\bar w})^{\rho-1},
\mbox{ with respect to }
  d \mu_\rho (z)  = \frac{\rho}{\pi} (1 +  |z|^2)^{-(\rho+1)}  dz.
\]
Note that the reference measure $\mu_\rho$ is typically {\it not} a
constant multiple of the invariant measure $\nu$. ${\check
K}_{\rho}$ is a projection kernel,  onto the orthonormal polynomials
$  \sqrt{{ \rho -1 \choose k}}  z^{k}$, $k = 0, \dots, \rho-1$ in
$L^2(\SS, \mu_\rho)$. In this case, $\rho$ is really the total
number of particles.

\vspace{.2cm}

\noindent {\bf Hyperbolic model.} Take $\Lambda = \UU$, the unit
disk, which we identify with the hyperbolic plane. Let
\[
   {\check K}_{\rho}(z, \bar w) =   \frac{1}{(1 - z \bar w)^{\rho+ 1}} \mbox{ with respect to  } d \mu_\rho(z) = \frac{\rho}{\pi} (1 - |z|^2)^{\rho - 1} dz,
\]
for any $\rho > 0$.  As in the planar model, $({\check K}_{\rho},
\mu_\rho)$ defines a determinantal process with infinitely many
points; the orthonormal polynomials in $L^2(\UU, \mu_\rho)$ being $
\sqrt{{\rho + k \choose k}} z^k$ for $k = 1,2, \dots$.

\vspace{.2cm}

The intensity measure in any determinantal processes is ${K}(z,z)
d \mu(z) $. In the above models we find that ${\check
K}_{\rho}(z,z)d  \mu_\rho(z)=\frac{\rho}{\pi} d \inv(z)$, where
$\inv$ is the {\bf invariant measure} on $\Lambda$, unique up to
constant. Here we use
 \be \label{meandens}
    d \inv_\CC (z) =  d z,  \ d\inv_\SS(z) =  \frac{d z}{ (1+|z|^2)^2} .  \  \ \mbox{ and }  d\inv_\UU(z)=\frac{\one_\UU\;  dz}{ (1-|z|^2)^2};
\ee The distribution of the above processes is invariant under
symmetries of the respective $\Lambda$, i.e.
linear fractional transformations of $\CC$ preserving the measure
$\inv_\Lambda$.

Our main theorem concerns the linear statistics for the point
process. Recall that for $f:\Lambda\to \R$, the ivariant measure
$\nu=\nu_\Lambda$ and the intrinsic gradient $\nablai$ we have
 $$
\|f\|_{H^1(\inv)}^2 =\int_\Lambda |\nablai f|^2 d\inv, \qquad
\|f\|_{L^1(\inv)} =\int_\Lambda |f| d\inv.
 $$
We say $f$ is in $H^1(\inv)$ or $L^1(\inv)$ if these corresponding
norm is finite.
 \begin{theorem}
\label{thm1}\label{main} For either the planar, spherical or
hyperbolic model, let $f\in H^1(\inv)\cap L^1(\inv)$.  Then, as
$\rho \ra \infty$, the distribution of
\[
 \sum_{z\in \mathcal Z} f(z)
 -  \frac{\rho}{\pi}\int_\Lambda f d\inv
\]
converges to a mean zero normal with variance $\frac{1}{4 \pi}
\|f\|_{H_1(\inv)}^2$.
\end{theorem}

Note that for both the limiting variance and the shift to make
sense it is necessary to have $f\in H^1\cap L^1$, so the theorem
holds for the most general test functions possible.

The fact that the variance is of order 1 manifests the advertised
repulsion.  The $H_1$-norm is conformally invariant, so one may
replace the intrinsic gradient and intrinsic measure by the planar
gradient and Lebesgue measure for the embedding.
%
%
%
%
%
%
%

This work is partially motivated by the recent results of Sodin and
Tsilerson  \cite{Sodin} on the three canonical {\em Gaussian
analytic functions} (GAFs) with zero sets invariant under the
symmetries of the plane, sphere, and hyperbolic plane. These
processes are also indexed by a density parameter, and \cite{Sodin}
establishes asymptotic normality for the corresponding linear
statistics, with $f \in C_0^2$. What is striking is that for GAFs
the variance actually decays as the density tends to infinity:
 $$ \Var \sum_{z\in \mathcal Z_\rho} f(z) \sim \rm{const}\times
\rho^{-1} \| \Delta_\iota f \|_{L^2(\inv)}^2.$$ Thus, the zeros of
typical GAFs are more orderly than their determinantal
counterparts.

The determinantal processes studied here, while attractive solely
on the basis of their invariance, also arise  as matrix models.
The planar case  is really just the ``infinite dimensional Ginibre
ensemble''. If $A$ is an $n \times n$ matrix of iid standard
complex Gaussians,  then as $n \uparrow \infty$  the point process
of $A$-eigenvalues converges to the planar model, and $\rho$  here
corresponds to scaling.  As for the spherical model, Krishnapur
\cite{Manju} has proved that it coincides with the eigenvalues of
$ A^{-1} B$, where $A$ and $B$ are independent $\rho \times \rho$
Ginibre matrices. Further, for integer $\rho$, Krishnapur provides
strong evidence that the hyperbolic points have the same law as
the singular points of $ A_0 + z A_1 + z^2 A_2 + \cdots $ in $|z|
< 1$ with again $A_0, A_1 , \dots$ independent $\rho \times \rho$
Ginibre matrices.

In all three cases, Krishnapur provides natural random analytic
functions for which $\mathcal Z$ is the set of zeros. Using an
integration by parts argument, Theorem \ref{main} can be
interpreted to say that the log absolute value
 of these analytic functions
converges to the Gaussian Free Field.  See Section 3 in
\cite{RidVir} for this relation in the Ginibre ensemble.

Theorem \ref{thm1} also identifies the present as a companion
paper to \cite{RidVir} which treats the limiting noise for the
Ginibre eigenvalues. The eigenvalues also define a determinantal
process in $\CC$, see \cite{Gin}. In \cite{RidVir} it is shown
that, along with an $H^1$-noise in the interior of $\UU$ similar
to above, there is an $H^{1/2}(\partial \UU)$ noise component in
the corresponding $n \uparrow \infty$ central limit theorem. The
invariance and lack of boundary effects in the three models
considered here makes for essentially different proofs that are
shorter and rely less on combinatorial constructions. The only
overlap is the cumulant formula which is the starting point for
both proofs.

The main Theorem \ref{main} is proved in three steps. In Section
\ref{genan} we establish some general conditions under which (smooth)
linear statistics are asymptotically normal, without computing the asymptotic variance. For
this, the fact that the kernels is an analytic projections and their
specific decay properties are crucial. In Section  \ref{ouran} we
check that these properties are satisfied by our models. In Section
\ref{morean} we determine the asymptotic variance and extend the
convergence to general test functions.

\section{General conditions for asymptotic normality}
\label{normal}\label{genan}

Taking a broader perspective, this section shows that under certain
conditions satisfied by the models we are considering, linear
functionals are asymptotically normal.

Let $B$ be a compact subset of $\CC$. Consider the following set-up.
 $K_\rho:B^2\to \reals$ is a set of kernels
indexed by $\rho$, which ranges in an unbounded subset of the
positive reals. The kernels here are the Hermitian and are {\bf
with respect to Lebesgue measure}; more precisely, if $\check
K_\rho$ denotes the kernels outlined above, and $\kappa=
d\mu_\rho(z)/dz$ is the density of the reference measure, then
here and below
 \be \label{Kdef}
 K_\rho(z,w)=\check K_\rho(z,w)
 \left(\kappa(z)\kappa(y)\right)^{1/2}.
 \ee


\subsection*{Kernel Properties}
\label{kprops}

The eventual asymptotic normality rests on the following asymptotic
properties of  $K_\rho$ as $\rho \to \infty$; all limit statements
and $o(\cdot)$ notations refer to this limit.  Throughout, $c$
denotes a numerical constants which may change from line to line.

\begin{itemize}
\item {\bf Uniform bound (UB).}   It holds
 \be
 \label{cond1}
 \|K_\rho\|_\infty:=\sup_{x,y\in B}|K_\rho(x,y)|\le c \rho.
 \ee
\item {\bf $L^1$ bound ($L^1$B).}
For $x,y\in B$  we have a bound \be |K_\rho(x,y)|\le
\varphi_\rho(x-y)  \mbox{  with  }
 \label{cond2}
 \left\|\varphi_\rho\right\|_1 \le c.
 \ee
 \item {\bf Interaction decay (ID).}
The above bounding function  satisfies
 \be
 \label{cond3}
 \||y|^{3} \varphi_\rho\|_1=o(\rho^{-1}).
 \ee
\item {\bf Limited local analytic projection  (LLAP) property.}
Assume that $B\subset \CC$. Fix $B_2$ compact so that $B_2\subset
B^0$. For $p=0,1,2$ we have \be
 \label{anpro}
\sup_{x,z\in B_2} \left| x^pK_\rho(x,z)  - \int_B K_\rho (x,y)\,y^p
\, K_\rho(y,z)\,dy \right|= o(1),
 \ee
Similarly, \be
 \label{canpro} \sup_{x,z\in B_2} \left| K_\rho(x,z)\bar
z^p  - \int_B K_\rho (x,y)\,\bar y^p \, K_\rho(y,z)\,dy \right|=
o(1).
 \ee
 \item{\bf Covariance (CO)}.
 For any function $F$ with bounded third derivatives and compact
support in the interior of $B$ we have $\Cov_\rho(
\partial_{z,\bar z}F, z\bar z) =o(1).$
\end{itemize}

Note of course that $\Cov_\rho(f, g)$ indicates the covariance of
$\sum_{z\in \mathcal Z} f(z)$ and $\sum_{z\in \mathcal Z} g(z)$  in
the $K_{\rho}$-process. The main proposition of this section (proved
as Proposition \ref{scum}) is:

\begin{proposition} \label{scum0}
Suppose that the above conditions are satisfied. Then, for the
corresponding determinantal process  any linear statistic $f$ with
compact support in the interior $B^o$ and bounded third derivatives
is asymptotically normal. We also have convergence of all moments.
\end{proposition}

That these conditions are satisfied by the planar, spherical and
hyperbolic models is delayed to the next section.  Here we provide a
lemma that sheds more light as to how condition LLAP arises.


\begin{lemma}[Analytic projections restrict to LLAP]
\label{LLAP} Let $\hat K_\rho:S^2\to \reals$ be a kernel for the
projection to the space of all analytic functions on the open set
$S\subset \CC$ with respect to measures $\mu_\rho$. For compact
$B\subset S$, let $K_\rho$ denote the restriction of $\hat
K_\rho(x,y) (\mu_{\rho}(x) \mu_{\rho}(y))^{1/2}$ to $B\times B$. If
$K_\rho$ satisfies UB and  ID, then $K_\rho$ satisfies the
LLAP.
\end{lemma}

\begin{proof}
Note that since for each $z$, the function $y\mapsto y K_\rho(y,z)$
is analytic, it follows from the analytic projection property that
$$
\int_S K_\rho(x,y)y^p K_\rho(y,z)\, dy = x^pK_\rho(x,z).
$$
Thus, for \re{anpro} it suffices to show that
 \be
 \label{outside}
 \sup_{x,z\in B_2}\left| \int_{S\setminus B} K_\rho(x,y)\,y^p \,K_\rho(y,z)\, dy  \right|= o (1),
 \ee
where recall $B_2 \subset B^o$. Setting $s=y-z$, there is a
polynomial $q$ of degree $p$ so that $|y^p| \le q(|s|)$ for all
choices of $z\in B_2, y\in B^c$. Also,    $q(|s|)\le c|s|^{3}$ as
soon as  $|s|$ bounded away from zero. So,  for $z$ a positive
distance from $S\setminus B$, we have that
$$
|y^p \,K_\rho(y,z)| \le c |s|^{3} \varphi_{\rho}(s),
$$
and by UB, ID, the absolute value of the left hand side of
$\re{outside}$ is bounded above by
$$
c \rho \int |s|^{3} \varphi_{\rho}(s)ds = c\times \rho \times
o(\rho^{-1}) =o(1).
$$
The proof of \re{canpro} is identical since $K_{\rho}$ is hermitian
symmetric.
\end{proof}
%
%

\subsection*{Cumulants}

Recall that for any random variable $X$, the cumulants $\Cum_k(X)$,
$k = 1, 2, \dots$, are the coefficients in the expansion of the
logarithmic generating function,
\[
   \log \ev [ e^{i t X} ] = \sum_{k=1}^{\infty}  \frac{(it)^k}{k!} \Cum_k(X),
\]
and $X$ is Gaussian if and only of $\Cum_k(X) = 0$ for all $ k \ge
3$. In any determinantal process $(K_{\rho}, \mu_{\rho})$, the
cumulants of the linear statistic $\sum f(z_k)$ have the explicit
form, \be \label{cumulants} \Cum_{k, \rho}(f) = \sum_{m=1}^m
\frac{(-1)^{m-1}}{m!}  \sum_{ k_1 + \cdots + k_m = k \atop k_1,
\dots, k_m \ge 1}
 \frac{k!}{k_1! \cdots k_m!}
                            \int   \Bigl[ \prod_{i=1}^m f(x_i)^{k_i} K_{\rho} (x_i, x_{i+1}) \Bigr]  dx_1 \dots dx_m,
\ee where  $x_{m+1} = x_1$ is understood, the integral ranges over
$m$ copies of the full space (here $B$), and again we are
absorbing the reference measure $\mu_{\rho}$ into the $K_{\rho}$
kernel. The structure behind formula (\ref{cumulants}) has been
employed in the past to establish asymptotic normality for
determinantal processes with various assumptions on the regularity
of $f$. See in particular the pioneering work of Costin-Lebowitz
\cite{CosLeb} and the later papers of Soshnikov, \cite{Sosh00} and
\cite{Sosh02}.  While going through cumulants, the method here is
quite different.

We define the multiple integrals: for $f$  a function of $x_1,\ldots
x_k$, \be \label{Int1} \KK_\rho^{\circ k}(f)=\int
\left[f(x_1,\ldots, x_k)\prod_{j=1}^k K_\rho(x_j,x_{j+1})
d\mu(x_j)\right] \ee (the indices are mod $k$), and, as another
shorthand, if the $f_i$ are all functions of one variable, we set
\be \label{Int2} \KK_\rho(f_1,\ldots, f_k)=\KK^{\circ
k}(f_1(x_1)\cdots f_k(x_k)). \ee Note that the cumulant
(\ref{cumulants}) is just a weighted sum of terms of the form
$\KK_\rho(g_1,\ldots, g_m)$, obtained by partitioning $\{1, 2,
\dots, k\}$ into $m$ parts $I_1, \dots, I_m$ of sizes $k_1, \dots
k_m$ and setting $g_i = \prod_{j \in I_i} f_j$.    Hence, we more
generally seek conditions for the vanishing of $\Cum_{\rho}(f_1,
\dots, f_k)$, defined in the obvious way,  for $k \ge 3$.

The first step is a collection of estimates on the integrals
(\ref{Int2}). The most fundamental of these are Lemmas 8 and 9
below.  The former allows one to reduce the dimension in certain
instances; the latter allows for the replacement of the test
functions $f_i$ by their cubic approximations.

\begin{lemma}\label{mainbound}
Assumptions $L^1$B and  UB imply that
$$
|\KK_\rho(f_1,\ldots,f_k)| \le c \rho \|f_1\|_1
\prod_{\ell=2}^\infty \|f_\ell\|_\infty.
$$
\end{lemma}

\begin{proof}
The integral is bounded above by
$$
   \int_{B^k}|f(x_1)|\,  \varphi_{\rho}(x_2-x_1)\cdots
   \varphi_{\rho}(x_k-x_{k-1})\,dx_1\cdots dx_k \times \|f_2\|_\infty\cdots \|f_\ell\|_\infty \times \sup_{x,y\in
   B}|K_\rho(x,y)|.
$$
Changing variables, $y_i=x_i-x_{i-1}$ for $i\ge 2$  allows the
remaining integral to be bounded  above by
  $$
  \int_{B\times (2B)^{k-1}}|f_1(x_1)|
  \varphi_{\rho}(y_2)\cdots \varphi_{\rho}(y_k)\,dx_1dy_2 \cdots dy_k=
  \|f_1\|_1\|\varphi_{\rho}  \|_1^{k-1},
  $$
and the claim follows from assumptions ($L^1$B, UB).
\end{proof}

\begin{lemma} \label{disjoint}
Assume $L^1$B, UB, ID and that $f_i$ are bounded for all $i$. If any
$f_i$ and $f_j$ are supported on disjoint compact sets,
then $\KK_\rho(f_1,\ldots, f_k)\to 0$.
\end{lemma}

\begin{proof}
Now use the bound
\begin{eqnarray*}
 \KK_\rho(f_1,\ldots, f_k) \le
   {c\|K\|_\infty
   \int \varphi_{\rho}(x_2-x_1)\cdots
   \varphi_{\rho}(x_k-x_{k-1})
  \,dx_1\cdots dx_k}
\end{eqnarray*}
where the integral on the right is over the product of the supports
of the $f_1$ to $f_k$.  By adjusting $c$, we may insert
$|x_i-x_j|^{3}$ (which is bounded below on the domain of
integration) to produce
\begin{eqnarray*}
 \KK_\rho(f_1,\ldots, f_k)
  & \le &c\|K_\rho\|_\infty
  \int_{B\times (2B)^{k-1} }
  \varphi_{\rho}(y_2)\cdots \varphi_{\rho}(y_k)|y_{i+1}+\ldots+y_j|^{3}
  \,dx_1dy_2 \cdots dy_k
  \\
&\le&
  c\|K_\rho\|_\infty \sum_{\ell={i+1}}^j
  \int_{B\times (2B)^{k-1} }
  \varphi_{\rho}(y_2)\cdots \varphi_{\rho}(y_k)|y_{\ell}|^{3}
  \,dx_1dy_2 \cdots dy_k
  \\
&\le& c \times \rho \times o(\rho^{-1}) = o(1), \qquad
\end{eqnarray*}
after a change of variables ($y_i = x_i - x_{i-1}, i > 1$) in line
one.
\end{proof}

Lemmas \ref{mainbound} and \ref{disjoint} were developed for the
following purpose.

\begin{lemma}\label{collapse}
Assume $L^1$B, UB, ID, and  LLAP.  Let  $f_1,\ldots, f_k$ be bounded
in $B$, with, for some $i$,  $f_i$ supported in a compact $B_1
\subset B^o$, and, for some $j \neq i$, $f_j(z)=z^p$ for $z\in B$,
with $p\in\{0,1,2\}$. Then
$$
\KK_\rho(f_1,\ldots, f_k) = \KK_\rho(f_1,\ldots,f_{j-1}\times
f_j,f_{j+1}\ldots, f_k)+o(1).
$$
Similarly, for $f_j(z)=\bar z^p$ we have
$$
\KK_\rho(f_1,\ldots, f_k) = \KK_\rho(f_1,\ldots,f_{j-1},f_j\times
f_{j+1}\ldots, f_k)+o(1).
$$
\end{lemma}

\begin{proof}
We prove the first claim with $f_j=z^p$;  the proof of the second
claim is identical. By the cyclic nature of $\KK_\rho$, we may
assume $i=1$. Fix a compact set  $B_2$ such that
 $B_1\subset B_2^o \subset B_2 \subset B^o $. By
the disjoint union decomposition
$$(B\times B)\setminus (B_2\times B_2) =\left(
(B\setminus B_2) \times B \,\right)\cup\left( \, B_2 \times
(B\setminus B_2)\right)$$ and Lemma \ref{disjoint} we have the
following (restrictions are placed only on functions with indices
adjacent to $j$):
\begin{eqnarray*}
   \lefteqn{
  \hspace{-5em}
\left|\KK_\rho(f_1,\ldots,
f_k)-\KK_\rho(f_1,\ldots,f_{j-1}\one_{B_2},f_j,f_{j+1}\one_{B_2},\ldots,
f_k) \right|     }\\
     &&\le \left|
     \KK_\rho(f_1,\ldots,f_{j-1}\one_{B\setminus B_2},f_j,f_{j+1},\ldots,
f_k)\right| \\&& \quad +\left|
     \KK_\rho(f_1,\ldots,f_{j-1}\one_{B_2},f_j,f_{j+1}\one_{B\setminus B_2},\ldots,
f_k)\right|  = o(1).
\end{eqnarray*}

Also,
\begin{eqnarray}
   \lefteqn{
  \hspace{-5em}
\KK_\rho(f_1,\ldots,f_{j-1}\one_{B_2},f_j,f_{j+1}\one_{B_2},\ldots,
f_k) - \KK_\rho(f_1,\ldots,f_{j-1}\one_{B_2}\times
f_j,f_{j+1}\one_{B_2},\ldots, f_k)}
 \label{llapu} \\ \label{sformula}
 &&=\int K_\rho(z_1,z_2)\cdots S(z_{j-1},z_{j+1})\cdots
K_\rho(z_k,z_1)\prod_{i=1\atop i\not=j}^k f_i(z_i)\, dz_i, \nonumber
\end{eqnarray}
where
$$
S(x,z)=\left(\int_B K_\rho(x,y)y^pK_\rho(y,z)\, dy -x^pK_\rho(x,z)
\right)\one(x,z\in B_2).
$$
But, by the LLAP assumption, we have that
$$
\sup_{x,z\in B} |S(x,z)| = o(1),
$$
and,
 after the familiar change of variables $y_i=x_i-x_{i-1}$  for $i\not=j+1$,
the argument used in Lemma \ref{mainbound} yields that the
difference \re{llapu} converges to zero. An identical application of
Lemma \ref{disjoint} gives
$$\KK_\rho(f_1,\ldots,f_{j-1}\one_{B_2}\times
f_j,f_{j+1}\one_{B_2},\ldots, f_k)-
\KK_\rho(f_1,\ldots,f_{j-1}\times f_j,f_{j+1},\ldots, f_k) = o(1),
$$
which concludes the proof.
\end{proof}

We close this subsection by showing that the assumed conditions
enable one to Taylor expand inside the $\KK_{\rho}$ integrals.

\begin{lemma}\label{taylor}
Assume that UB, $L^1$B and ID hold. Let $f_j$, $1\le j \le k$ have
bounded third derivatives. Then we have
$$
\KK_\rho(f_1,\ldots,f_k)=\sum_{m=0}^2
\sum_{m_2+\ldots+m_k=m}\KK^{\circ
k}_\rho\left[f_1(x_1)\bigotimes_{i=2}^k \left(f_i^{(m_i)}(x_1)
(x_i-x_1)^{\otimes m_i}\right)\right]+o(1) ,$$ where we use the
standard tensor notation for the full first and second derivatives.
\end{lemma}

\begin{proof}
Starting at $\ell=2$, we will step-by-step replace
$f_1(x_1)f_2(x_2)\cdots f_\ell(x_\ell)$, $\ell\ge 2$ by an
approximation $g_\ell$ of degree 2 at $x_1$:
$$g_\ell(x_1;x_2,\ldots,x_\ell) = f_1(x_1) \sum_{m=0}^2
\sum_{m_2+\ldots+m_\ell=m}\left[\bigotimes_{i=2}^\ell
\left(f_i^{(m_i)}(x_1) (x_i-x_1)^{\otimes m_i}\right)\right]
$$
Note that all the $g_\ell$ are bounded on $B^{\ell}$. For the
step-by-step replacement procedure we  need to bound
$d_\ell=g_{\ell-1} f_\ell- g_{\ell}$. Towards this end, let
$$
f_\ell^*=f_\ell^*(x_1,x_\ell)=f_\ell(x_1)+f_\ell'(x_1)(x_\ell-x_1)+f_\ell''(x_1)(x_\ell-x_1)^{\otimes
2}.
$$
Certainly,
$$|f_\ell^*(x_1,x_\ell)-f_\ell(x_\ell)| \le c |x_1-x_\ell|^3,
$$
for $C$ independent of $\ell,x_1,x_\ell$. Also,
$$|d_\ell|=|g_{\ell-1} f_\ell- g_{\ell}| \le
|g_{\ell-1}| |f_\ell-f^*_\ell|+|g_{\ell-1} f^*_\ell-g_{\ell}|.
$$
Let $y_\ell=x_\ell-x_1$. Since $g_{\ell-1}$ is bounded,  we have
 \be
 \label{firstbound}
\left|g_{\ell-1}||f_\ell-f^*_\ell\right|\le c \left(|y_2|^{3}+
\ldots +|y_{\ell}|^{3}\right)
 \ee
for the range of $y_i$. As $g_{\ell}$ is produced from  $g_{\ell-1}
f^*_{\ell}$ by dropping all terms that are of degree 3 or 4 in the
$y_i$, there is a constant $c$ such that
 \be
 \label{secondbound}
 |g_{\ell+1}-g_\ell f^*_{\ell+1}|
 \le c
 \left(|y_2|^{3}+\ldots+|y_{\ell}|^{3}\right)
 \ee
 on the range of the $y_i$.
(Any monomial in $y_i$ of degree 3 or 4 and coefficient 1 is bounded
above on the compact range by the right hand side for $c$ large
enough).

Now we write the difference
\begin{eqnarray}
   \lefteqn{
  \hspace{-5em}
   \big|\KK_\rho
   \left[g_{\ell-1}(x_1;x_2,\ldots ,x_{\ell-1})
   f_{\ell}(x_\ell)\cdots f_k(x_k)
   \right]}\nonumber
     \\
     &&\qquad -\;\KK_\rho\left[ g_{\ell}(x_1;x_2,\ldots,x_{\ell})
     f_{\ell+1}(x_{\ell+1})\cdots f_k(x_k)\right]\big| \nonumber \\
   &= & \left|\KK_\rho\left[ d_\ell(x_1;\,x_2,\ldots,x_{\ell})
   f_{\ell+1}(x_{\ell+1}) \cdots f_k(x_k)
   \right]\right|.  \label{kkdif}
\end{eqnarray}
By UB and $L^1$B \re{kkdif} is bounded above by
\begin{eqnarray*}
   \lefteqn{c\|K_\rho\|_\infty \int_{B^{k}}
  \varphi_\rho(x_2-x_1)\cdots\varphi_\rho(x_{k}-x_{k-1})|d_\ell(x_1;\,x_2,\ldots,x_{\ell+1})|
  dx_1\cdots dx_k}
  \\ & \le & c\|K_\rho\|_\infty
  \int_{B\times (2B)^{k-1}}
  \varphi_\rho (y_2)\cdots\varphi_{\rho}(y_k)  \Bigl( \sum_{i=2}^{\ell} |y_i|^3 \Bigr)
  dx_1dy_2 \cdots dy_k,
\end{eqnarray*}
with \re{firstbound} and \re{secondbound} used in the second line.
Again by $L^1$B and ID, this in turn is upper bounded by
$$
  c |B| \|K_\rho\|_\infty \left\|\,
  \varphi_\rho \times|y|^{3}\,
  \right\|_1 \le c \rho \times o(\rho^{-1})=o(1)
$$
as required.
\end{proof}

%
\subsection*{Proof of the proposition}

The above bounds on $\KK_{\rho}$ made use of UB, $L^1$B , ID, and
LLAP.  If we add CO to the mix, the result is the following.

\begin{proposition}\label{scum}
Assume that $K_\rho$ satisfies conditions UB, $L^1$B, and  ID. For
$k \ge 3$, let $f_1,\ldots f_k$, be of compact support and have
continuous third derivatives.  If in addition CO holds, then
$\Cum_\rho(f_1,\ldots,f_k) \to 0$.
\end{proposition}


To prove this, remember that $\Cum_\rho(f_1,\ldots,f_k) $ is a
weighted sum of terms in the form $ \KK_\rho(g_1,\ldots,g_m)$ each
$g_j$ being a product of the underlying $f$'s (here $m \le k$).
 Lemma \ref{taylor} gives
 \be
 \label{kk}
\KK_\rho(g_1,\ldots,g_m)=\sum_{\ell=0}^2
\sum_{\ell_2+\ldots+\ell_k=m}\KK^{\circ
m}_\rho\left[g_1(z_1)\bigotimes_{i=2}^k \left(g_i^{(z_i)}(\ell_1)
s_i^{\otimes \ell_i}\right)\right]+o(1),
 \ee
 in which $s_i=z_i-z_1$, and we use  complex coordinates $s_i,\bar
s_i$. For example, $$g_i^{(2)}(z_1) s_i^{\otimes 2}=
\partial_z\partial_z g_i (z_1)s_i^2+
\partial_{\bar z}\partial_{\bar z} g_i (z_1)\bar s_i^2+2
\partial_z \partial_{\bar z} g_i (z_1) s_i \bar s_i.$$
Further, the $k$-fold integrals on the right hand side of \re{kk}
are all of the  form $\tilde K^{\circ k}_\rho (h(z_1)\sigma_i
\sigma_j)$ where $h(z_1)$ is a $C^3$ compactly supported (in $B$)
function of $z_1$ and $\sigma_\eta=s_\eta,\bar s_\eta$ or $1$. Since
functions of at most three of the $z_i$-s are present in this
integrand, Lemma \ref{collapse} with $p=0$ allows us to reduce it to
an at most threefold integral:
$$\KK^{\circ k}_\rho (h(z_1)\sigma_i
\sigma_j)=\KK_\rho^{\circ (d+1)} (h(z_1)\sigma_i \sigma_j)+o(1),$$
where $d = 0, 1$ or $2$  is the number of distinct variables in
$\sigma_i \sigma_j$.

 Next,  two applications of Lemma \ref{collapse} gives
\begin{eqnarray*}
\KK_{\rho}^{\circ 2}(h(z_1)s_2)&=&
\KK_{\rho}(h,z)-\KK_{\rho}(hz,1) \\
&=&\KK_{\rho}(h,z)-\KK_{\rho}(hz)+o(1) =o(1).
\end{eqnarray*}
Similarly, all except the $s_i\bar s_j$ terms vanish. Even among
those, only half of the terms with $i\not=j$ survive, depending on
the order of conjugation. Again by successive applications of Lemma
\ref{collapse}
$$
\KK_{\rho}^{\circ 3}(h(z_1)s_2\bar s_3)=\KK_{\rho}^{\circ 2}(\bar
s_1 h(z_1)s_2)+o(1)=o(1),
$$
while \be \label{remains} \KK_{\rho}^{\circ 3}(h(z_1)\bar s_2
s_3)=\KK_{\rho}^{\circ 3}(h(z_1)\bar s_2 s_2) + o(1) =
\KK_\rho(h,z\bar z)-\KK_\rho(h z \bar z) + o(1). \ee Apart from the
error,  the latter equals $- \Cov_{\rho}(h, z \bar z)$.

Therefore, aside from a constant term,  the only $O(1)$
contributions to the cumulant sum are of the form (\ref{remains}).
The possible choices of $h$ are: with $F = f_1 f_2 \cdots f_k$,
\begin{eqnarray*}
G_i&=&F\times (\partial_{z,\bar z} g_i)/g_i\\
G_{ij}&=&F\times (\partial_{z} g_i)(\partial_{\bar z} g_j)/(g_ig_j),
\end{eqnarray*}
and our full $\KK_{\rho}$ formula reduces to $\KK_\rho(\cdot,z\bar
z)-\KK_\rho(\cdot \times z \bar z)$ applied to
$$
\sum_{i=2}^m G_i + \sum_{2\le i<j\le m} G_{ij}.
$$
Reverting back to the original test functions $f_1, f_2, \dots,
f_k$, this is a weighted sum of the functions
\begin{eqnarray*}
F_u&=&F\times (\partial_{z,\bar z} f_u)/f_u\\
F_{uv}&=&F\times (\partial_{z} f_u)(\partial_{\bar z} f_v)/(f_uf_v).
\end{eqnarray*}
Since $\Cum_{\rho}(f_1, \dots, f_k)$
 is symmetric under  permutations of
indices of the $f_i$'s,  it suffices to show that the total weight
over the cumulant sum for each one of the two types of terms $F_u$
and $F_{uv}$ vanishes.
Also, as each $g_i$ is a product of $k_i$ of the $f_i$, then $G_i$
is a sum of $k_i$ terms of type $F_{u}$ and $k_i(k_i-1)$ terms of
type $F_{uv}$. Similarly, when $i\not=j$, $G_{ij}$ is a sum of
$k_ik_j$ terms of type $F_{uv}$. Thus, the total number of terms of
each type is given by
\begin{eqnarray*}
\mbox{type }F_u&:& k_2+\ldots +k_m\\
\mbox{type  }F_{uv}&:&\sum_{i=2}^m k_i(k_i-1) + \sum_{2\le i<j\le
m}k_ik_j,
\end{eqnarray*}
while each type appears in the cumulant sum with a different
coefficient.

Finally we invoke property CO: $\KK_\rho(\partial_{z,\bar z}F,z\bar
z)-\KK_\rho((\partial_{z,\bar z}F) \times z \bar z)\to 0 $.  Since
 $$ \partial_{z,\bar z}F  =  F \times \Bigl( \sum (\partial_{z,\bar z} f_u)/f_u
  + \sum
 \partial_{z} f_u(z) \partial_{\bar z} f_v / (f_u f_v) \Bigr),
 $$ after  subtracting $\partial_{z,\bar
z}F/k$ from each term of type $F_u$ it  can be replaced by $k-1$
terms of type $F_{uv}$ with the  opposite sign. Thus,  our final
count is
\begin{eqnarray}
\label{finalcount} \mbox{type }F_{uv}&:& \sum_{i=2}^m (k_i^2-kk_i) +
\sum_{2\le i<j\le m}k_ik_j.
\end{eqnarray}
That is to say, each $m \le k$ term in the cumulant sum is the same
constant multiple of (\ref{finalcount}). That this vanishes for
$k\ge 3$ when summed over the full cumulant expansion is the content
of the next lemma.

\begin{lemma}
\label{vanish} For each $m\ge 1$ let $\varphi(k,m,k_1,\ldots ,k_m)$
be a real-valued function. With
$$
\Upsilon_k(\varphi) =\sum_{m=1}^{k} \frac{(-1)^{m-1}}{m}
\sum_{k_1+\ldots+k_m=k\atop k_1,\ldots ,k_m\ge
1}\frac{\varphi(k,m,k_1,\ldots,k_m)}{k_1!\cdots k_m!},
$$
it holds that
\begin{eqnarray}
\Upsilon_k\left(\sum_{i=2}^k(k_i^2 -kk_i)+\sum_{2\le i<j\le m}k_ik_j
\right)=0. \label{gff}
\end{eqnarray}
for all $k \ge 3.$
\end{lemma}

\begin{proof}
First realize, if we denote $\varphi\equiv 1 $ by 1, then
\begin{eqnarray*}
\log(e^x)&=&\sum_{m=1}^\infty \frac{(-1)^{m-1}}{m}\left(
\frac{x}{1!}+\frac{x^2}{2!}+\frac{x^3}{3!}+\ldots\right)^m\\
&=& \sum_{k=1}^\infty \Upsilon_k(1)\,x^k,
\end{eqnarray*}
which explains why the $k \ge 3$ cumulant of any constant is zero.
Now set
\begin{eqnarray*}
f=f(x,y)= e^{x}+xye^{x}&=&  1+
\frac{x(1+y)}{1!}+\frac{x^2(1+2y)}{2!}+\frac{x^3(1+3y)}{3!}+\ldots
\end{eqnarray*}
so that the coefficient of the $y$ term in the $y$-power series
expansion of $\log f$ is
\begin{eqnarray*}
\left.\frac{d}{dy}(\log f)\right\vert_{y=0}=\sum_{k=1}^\infty
\Upsilon_k(k_1+\ldots+k_m) \,x^k
\end{eqnarray*}
and similarly,
\begin{eqnarray*}
\frac{1}{2}\left.\frac{d^2}{d^2y}(\log f)\right\vert_{y=0}=
\sum_{k=1}^\infty \Upsilon_k\left(\mbox{$\sum_{1\le i<j\le
m}k_ik_j$}\right)\,x^k.
\end{eqnarray*}
In order to obtain the pure quadratic sums  we set
\begin{eqnarray*}
g=g(x,y)=e^{x}+y(xe^x+x^2e^{x})&=&1+\frac{x(1+y)}{1!}+\frac{x^2(1+2^2y)}{2!}+\frac{x^3(1+3^2y)}{3!}+\ldots.
\end{eqnarray*}
The coefficient of the $y$ term in this  power series expansion
reads
\begin{eqnarray*}
\left.\frac{d}{dy}(\log g)\right\vert_{y=0}=\sum_{k=1}^\infty
\Upsilon_k(k_1^2+\ldots+ k_m^2)\,x^k.
\end{eqnarray*}

These series produce the types of terms we are after up to the fact
that  our cumulant expressions do not have the first coeffiicient
$k_1$. To omit this, the above may be modified as in
\begin{eqnarray*}
s_1&=&\left.\frac{d}{dy}\frac{(\log f)(e^x-1)}{f-1}\right\vert_{y=0}
= \left.\frac{d}{dy} \sum_{m=1}^\infty
\frac{(-1)^{m-1}}{m} (e^x-1)(f-1)^{m-1}\right\vert_{y=0}\\
&=& \sum_{k=1}^\infty \Upsilon_k(k_2+\ldots+k_m)\,x^k,
\end{eqnarray*}
and
\begin{eqnarray*}
s_2=\left.\frac{d}{dy}\frac{(\log
g)(e^x-1)}{g-1}\right\vert_{y=0}=\sum_{k=1}^\infty\Upsilon_k(k_2^2+\ldots+
k_m^2)\,x^k,
\end{eqnarray*}
and finally
\begin{eqnarray*}
s_{11}=\frac{1}{2}\left.\frac{d^2}{d^2y}\frac{(\log
f)(e^x-1)}{f-1}\right\vert_{y=0} = \sum_{k=1}^\infty
\Upsilon_k\left(\mbox{$\sum_{2\le i<j\le m}k_ik_j$ }\right)\,x^k.
\end{eqnarray*}
Now we have an easy proof of the claim. Simply note that the right
hand side of \re{gff} equals the coefficient of the $x^k$ term in
$$
-x\, \frac{d}{dx}s_1+s_2+s_{11} =\frac{ x^2}{2};
$$
the latter being a straightforward computation.
\end{proof}

\section{Properties satisfied by our models}
\label{ouran}

We now verify that the conditions for asymptotic normality of
Section \ref{kprops} hold for the three invariant models.  First
the simple bounds UB, $L^{1}B$ and ID are checked.  Note that for
the planar and hyperbolic models, property LLAP follows from Lemma
\ref{LLAP}: both families of kernels are projections onto the
space of analytic functions on the plane or unit disk. For the
spherical model, in which the kernel projects onto a set of finite
degree, LLAP requires separate proof.

%

The first three properties follow from a similar pointwise bound for
each of the kernels in question. Recall the definition \re{Kdef} of
$K_\rho(z,w)$.

\bigskip

\noindent {\bf Planar model.} We have
$$
K_\rho(z,w)=\frac{\rho}{\pi}e^{\,\rho( z\bar w - |z|^2/2+|w|^2/2)},
$$
and so \be |K_\rho(z,w)|=\frac{\rho}{\pi}e^{\,\rho(\Re(z\bar w) -
(|z|^2+|w|^2)/2)} = \frac{\rho}{\pi} e^{-\rho|z-w|^2/2}. \label{up1}
\ee That is,  $\varphi_\rho(z)=\frac{\rho}{\phi} e^{-\rho|z|^2/2}$,
and it is immediate  that conditions UB, $L^1$B,  and ID are
satisfied for all $z,w\in \CC$.

\bigskip
\noindent {\bf Hyperbolic model.} Denote $t=(1-|z|^2)(1-|w|^2)$ and
$s=(1-z\bar w)^2$. Then $|s|=t+|z-w|^2$, and we have
 $$
K_\rho(z,w)= \frac{\rho}{\pi s}(t/s)^{(\rho-1)/2}.
 $$
 If  $|z|<R<1$, then  $|s|\in[(1-R)^2,4]$, and, assuming
$\rho\ge 2$,  we get
\begin{eqnarray}
\label{up2} |K_\rho(z,w)|&=&\frac{\rho}{\pi |s|}
|t/s|^{\frac{\rho-1}{2}} = \frac{\rho}{\pi |s|}
\left[1-\frac{|z-w|^2}{|s|}\right]^{\frac{\rho-1}{2}}
\\&\le&
\frac{\rho}{(1-R)^2}\left[1-\frac{|z-w|^2}{4}\right]^{\frac{\rho-1}{2}}\le
\frac{\rho}{(1-R)^2}\;e^{-(\rho-1)|z-w|^2/8}=:\varphi_\rho(z-w).
\nonumber
\end{eqnarray}
Thus,  conditions UB, $L^1$B, and ID are satisfied for all $|z|\le
R$, $|w|<1$.

\bigskip
\noindent {\bf Spherical model.} Now put $t=(1+|z|^2)(1+|w|^2)$,
$s=(1+z\bar w)^2$. Then $t=|s|+|z-w|^2$, and the kernel reads
$$
K_\rho(z,w)= \frac{\rho}{\pi t}(s/t)^{(\rho-1)/2} .$$ If,  both
$|z|$ and $|w|<R$, then $t\in[1,b]$ with $b=(1+R^2)^2$, and
\begin{eqnarray}
\label{up3} |K_\rho(z,w)|&=&\frac{\rho}{\pi t}
|s/t|^{\frac{\rho-1}{2}} = \frac{\rho}{\pi t}
\left[1-\frac{|z-w|^2}{t}\right]^{\frac{\rho-1}{2}}
\\&\le&
\frac{\rho}{\pi t}\;e^{-(\rho-1)|z-w|^2/(2t)}\le
\rho\;e^{-(\rho-1)|z-w|^2/(2b)}=:\varphi_\rho(z-w) \nonumber
\end{eqnarray}
as soon as $\rho \ge 2$.  Just as before, we have UB, $L^1$B, and
ID when $|z|,|w|\le R$.

Rounding out the basic properties we have:

\begin{lemma} The spherical model restricted to a ball $\{z\le |B|\}$ has the LLAP property.
\end{lemma}

\begin{proof} Fix $p>0$, and assume that $\rho \ge 1+p$.
We consider a truncated kernel, which is a projection to the space
of polynomials of degree at most $\rho-1-p$ with respect to the same
measure as $K_\rho$, that is $\mu_\rho$. That is, we introduce
$$
\frac{\pi t^\frac{\rho+1}{2}}{\rho} \hat
K_\rho(z,w)=\sum_{k=0}^{\rho-1-p}{\rho-1 \choose k}(z\bar w)^{k},
$$
with again $t =(1-|z|^2)(1-|w|^2)$.  This truncated kernel is shown
to have LLAP, and then the truncation is shown to make a negligible
difference.

First note that
 \be
 \label{sl1}
 \int_\CC K_\rho(x,y)y^p\hat K_\rho(y,z) dy = x^p\hat K_\rho(x,z),
 \ee
since $K_\rho$ is a projection into the space of polynomials of degree at
most $\rho-1$ and for $z$ fixed and $y^p\hat K_\rho(y,z)$ is
$\mu_\rho(y)^{1/2}$ times such a polynomial. Also,  since for $z$
fixed, $K_\rho(y,z)-\hat K_\rho(y,z)$ is orthogonal to all powers at least
$\rho-1-p$ of $\bar y$ on any radially symmetric set, we
 have
 \be
 \label{sl3}
 \int_{|y|\le B} K_\rho(x,y)y^p( K_\rho(y,z)-\hat
K_\rho(y,z)) dy =0.
 \ee

Next, using the bound ${\rho-1 \choose k} \le {\rho-1-p \choose
k}\rho^p$ we find
$$\frac{\pi t^\frac{\rho+1}{2}}{\rho} |\hat K_\rho(z,w)|\le \rho^p(1+|zw|)^{\rho-1-p}.$$
And, when $|z|<b<B<|w|$, we have  $(1+|zw|)^2/t<a^2$ for some $a<1$.
Thus,
 $$|\hat K_\rho(z,w)|
 \le \frac{\rho^{p+1}}{t^{p/2+1}}\left[\frac{1+|zw|}{t^{1/2}}\right]^{\rho-1-p}\le\frac{\rho^{p+1}}{t^{p/2+1}} a^{\rho-1-p}.$$
Similarly, for $|z|<b<|B|<|w|$, we have
$$|K_\rho(z,w)|\le \frac{\rho}{t}\left[\frac{1+|zw|}{t^{1/2}}\right]^{\rho-1}\le\frac{\rho}{t} a^{\rho-1}.$$
By Hermitian symmetry this implies that
 \be
 \label{sl2}
 \sup_{|x|,|z|<b}\left|\int_{|B|<|y|} K_\rho(x,y)y^p
\hat K_\rho(y,z) dy\right| \le \rho^{p+2} a^{2\rho-2-p}
\int_{|B|<|y|}\frac{|y|^p}{(1+|y|^2)^{p/2+2}}\, dy= o(1).
 \ee
Finally, an easy estimate shows
 \be
 \label{sl4}
 \sup_{|x|,|z|<b} |x^p K_\rho(x,z)-x^p\hat K_\rho(x,z)|= o(1),
 \ee
 and  \re{sl1},\re{sl2},\re{sl3},\re{sl4} together imply the first
part of the LLAP property.  The second part follows from Hermitian
symmetry.
\end{proof}

\section{Asymptotic variance and general test functions}
\label{morean}

Our goal is to prove asymptotic normality with explicit variances
for any $f\in L^1\cap H^1$. We do this by proving normality and
determining the variance asymptotics for an $\|\cdot \|_{H^1}$-dense
set of functions and then giving a uniform variance bound for all
functions.

First note the general formula valid for all bounded  $f,g$  of
compact support:
\begin{eqnarray*}
  \Cov_{\rho}(f, g) & = & \int f(z) g(z) K_{\rho} (z,z) \, dz   - \int  \int f(z) g(w) K_{\rho} (z,w) K_{\rho} (w,z) \, dz dw \\
                                & = &  \int   \int  f(z) (g(z) - g(w)) | K_{\rho}(z, w)|^2 \, dz
                                dw,
\end{eqnarray*}
which, after symmetrization, reads \be \label{cov}
 \Cov_{\rho}(f, g)
                                =    \frac{1}{2} \int \int  (f(z) - f(w)) (g(z) - g(w)) | K_{\rho}(z, w)|^2 \, dz
                                dw.
\ee

\begin{lemma}
The subset of smooth functions with compact support not containing
$\infty$ is  $\|\cdot \|_{H^1}$-dense  in $H^1(\nu)\cap L^1(\nu)$.
\end{lemma}

Note that this subset is {\it not} dense in  $H^1(\UU)$, only among
$H^1\cap L^1$ functions: harmonic functions $h$ in $\UU$ are
$H^1$-orthogonal to any compactly supported $f$. This can be seen
via an integration by parts, moving the gradient from $f$ to produce a
$\triangle h = 0$.

\begin{proof}

 Replacing $f$ by $(f\wedge  b)\vee(-b)$ and letting $b\to
\infty$ shows that bounded functions are dense. Then convolving a
bounded $f$ with a smooth probability density approaching $\delta_0$
shows that bounded $C^3$ functions are dense.

First consider the planar or hyperbolic case, and use the
invariant gradient $\nablai$, measure $\nu$ and distance $\disti$.
We may apply a sequence of smooth cutoff functions $g_r$ to $f$
which are equal to $1$ on the ball of radius $r$ but are compactly
supported and have $|\nablai g_r|\le 1$. Let $h_r=1-g_r$. We have
$$\nablai (f - fg_r) = h_r \nablai f + f \nablai h_r
$$
and therefore
\begin{eqnarray*}
\|f-f g_r\|^2_{H^1} &\le &2 \int |h_r\nablai f|^2+ |f \nablai
h_r|^2\,d\nu(z)  \\
&\le& 2\int_{{\rm \dist}_\iota(0,z)>r} |\nablai f(z)|^2\, d\nu(z) +
2 \|f\|_\infty \int_{{\rm \dist}_\iota(0,z)>r} |f(z)|\,d\inv(z)
\end{eqnarray*}
these  converge to $0$ for bounded $f\in L^1\cap H^1$ as $r\to
\infty$.

For the sphere, we again consider smooth $f$, and note that adding a
constant to $f$ does not change its $H^1$-norm (formally, the space
$H^1$ consists of equivalence classes of functions which differ by a
constant). Adding a constant does not change the fact that $f\in
L_1$ either, as the invariant measure $\nu_\SS$ is finite. So we may
assume that $f(\infty)=0$, and by smoothness and compactness
$f(z)\le c_f\,\disti(\infty,z)$. We now take
$g_\eps(z)=((\disti(z,\infty)/\eps-1)\vee 0 )\wedge 1$, which is
supported on points at least $\eps$ away from $\infty$. Also,
$|\nablai g(z)|\le 1/\eps$ and vanishes for $z$ more than $2\eps$
away from $\infty$. As before, we have
\begin{eqnarray*}
\|f-f g_\eps\|^2_{H^1} &\le & 2\int_{{\rm
dist}_\iota(\infty,z)<2\eps} |\nablai f(z)|^2\, d\nu(z) +
2\int_{{\rm dist}_\iota(\infty,z)<2\eps} \frac{c_f^2
\,\disti(\infty,z)^2}{\eps^2} \,d\nu(z)
\end{eqnarray*}
Both terms converge to 0 when $\eps\to 0$, as required.
\end{proof}

\begin{lemma}[Asymptotic variance for a dense set]
\label{conc}
 Let $f$ and $g$ be $C^1$ and of compact support in $\Lambda$, where
$\Lambda = \CC$ for the plane or the sphere or $\UU$ for the
hyperbolic plane.  Then \begin{eqnarray}
\label{cov1}
   \Cov_{\rho}(f,g)
                               & \ra   &  \frac{1}{4 \pi} \langle f,g
                               \rangle_{H_1}.
\end{eqnarray}
\end{lemma}

\begin{proof}
It suffices to compute $\lim_{\rho \ra \infty} \Var_{\rho}(f)$ for
$f \in C_0^1$ as the covariance may be identified by substituting
$f+g$ for $f$.  First, by Taylor's theorem with remainder there is a
bounded non-negative function $\ep(r)$ tending to $0$ as $r
\downarrow 0$ for which
\begin{eqnarray*}
 \lefteqn{ \hspace{-1.5cm}
   \Bigl|  \Var_{\rho}(f)  -  \int_B \int_B  ( \nabla f(z) \cdot (w-z) )^2  | K_{\rho}(z, w) |^2 \, dz dw \Bigr|  } \\
    & \le & \int_B \int_B \ep(|z-w|)  |z-w|^2 \phi_{\rho}^2(z-w) \, dz dw = o(1).
\end{eqnarray*}
Now examine the remaining integrand
\[
    I(z, w) :=  | \nabla f(z)|^2 |z-w|^2 \cos^2(\theta)  | K_{\rho} (z, w) |^2,
\]
where $\theta(z,w)$ is the angle between $f(z)$ and $w-z$, under the
change of variables $w = z + \rho^{-1/2} w'$.  Pointwise, in each of
the three models, we have
\[
  \frac{1}{\rho^2}  | K_{\rho}(z, z + \frac{w'}{\sqrt{\rho}} ) |^2 \ra
   \frac{1}{\pi^2 \psi(z)^2} \exp \Bigl( - |w'|^2/ \psi(z) \Bigr),
\]
where $\psi(z) = 1$ (plane), $= 1 + |z|^2$ (sphere), $ = 1 - |z|^2$
(hyperbolic plane). This would result in the limiting formula for
the variance: with $\theta'$ denoting the limiting angle between $z$
and $w'$,
\begin{eqnarray}
\label{lastthing} \lefteqn{ \Var_{\rho}(f) \ra
  \frac{1}{2\pi^2} \int_{B}   |\nabla f(z)|^2  \int_{\CC} \cos^2(\theta')   \frac{|w'|^2}{\psi(z)^2} e^{- |w'|^2/ \psi(z)}  \, dw' dz} \\
  & = & \frac{1}{2 \pi^2}  \int_{B}   |\nabla f(z)|^2  \int_{0}^{\infty}  \int_0^{2\pi} \cos^2(\theta')  {|w'|^3} e^{- |w'|^2} \, d \theta'
    d |w'| d z
    =    \frac{1}{4 \pi}  \int_{B}   |\nabla f(z)|^2  \, dz. \nonumber
\end{eqnarray}
On the other hand, $z$ and $w$ ranging in a bounded set and
$||\nabla f||_{L^{\infty}} < \infty$, the right hand side of
\[
   | I(z, z + \frac{w'}{\sqrt{\rho}}) |   \le c |w'|  e^{ - |w'|^2/c}
\]
is integrable on $B \times \CC$ (this again uses (\ref{up1}),
(\ref{up2}), and (\ref{up3})).  Therefore, dominated convergence
validates (\ref{lastthing}) and completes the proof.
\end{proof}

\begin{corollary}\label{coco} The property CO holds for out models: $\Cov_{\rho}(\partial{z\bar z} f, z \bar z \one_{B}(z) ) \to 0$
for $f\in C^3$, ${\supp(f)}\subset B^o$ and $B$ compact.
\end{corollary}

\begin{proof}
Consider a smooth $g$ satisfying $\one_{B'}\le g\le \one_B$, where
$B'$ is a neighborhood of  $\supp(f)$. Lemma $\re{conc}$ with
implies
 \begin{eqnarray*} {4 \pi}\,
\Cov_\rho(\partial_{z\bar z} f, g\, z\bar z)\to \langle
\partial_{z\bar z} f, g z\bar z\rangle_{H_1} = - \langle f,
\partial_{z\bar z} (g z\bar z)\rangle_{H_1}=0.
 \end{eqnarray*}
For the last equality, note that $\nabla f$ and $\nabla
\partial_{z\bar z} (g z\bar z)$ have disjoint support.

 Now
$$\Cov(\partial_{z\bar
z} f,(\one_B-g))=\KK_\rho(\partial_{z\bar z}f,
(\one_B-g))-\KK_\rho(\partial_{z\bar z}f\times (\one_B-g)),$$ where
the second term vanishes, and the first one converges to $0$ by
Lemma \ref{disjoint} and the fact that the arguments have disjoint
support.
\end{proof}

\begin{lemma}[Variance bound]
\label{concbound}
 Let $f\in H^1(\nu)\cap L^1(\nu)$. There is a universal $c>0$ so that for all $\rho>1$ we
have \be \label{cov2}
  \Var_{\rho}(f)
                                \le      c \int_{\Lambda}  |\nablai f(z)|^2  \, d\nu(z).
\ee
\end{lemma}

\begin{proof}
By considering the negative and positive parts of $f$ separately, we
may assume $f\ge 0$.
In each of the three models by the invariant version of  \re{cov}
we have, for $f$ bounded and compact support
 \be\label{varfor} \Var_{\rho}(f) = \frac{1}{2}
\int_{\Lambda} \int_{\Lambda}  |f(z) - f(w) |^2 |\Ki(z,w)|^2
d\nu(z)\, d\nu(w),
 \ee
where $\Ki(z,w)=K(z,w)(\eta(z)\eta(w))^{-1/2}$ and
$\eta=d\nu(z)/dz$ is the density of the invariant measure
\re{meandens}.

Repeating the identities in \re{cov} for the invariant $\Ki$, and
using that for $\rho$ fixed $\Ki$ is bounded, we get that
\re{varfor} extends to all $f\in L^2(\nu)$, in particular for
bounded $f \in L^1(\nu)$.  Now replace the nonnegative $f\in
L^1(\nu)$ by $f_n=f\wedge n$. Let $\mathcal V(f)$ denote the right
hand side of \re{varfor}.  Since $|f_n(z)-f_n(w)|$ is monotone
increasing in $n$, the monotone convergence theorem gives
$\mathcal V(f_n)\to \mathcal V(f)$. We also have $\Var_\rho(f_n)
\to \Var_\rho(f)$: the mean converges to a finite limit as $f\in
L^1(\nu)$ and the second moment converges by the monotone
convergence theorem. Thus $\re{varfor}$ holds for all $f \ge 0$ in
$L^1(\nu)$, although a priori both sides may be infinite.

For isometries $T$ of $\Lambda$ (i.e. linear fractional
transformations preserving  $\nu$) we have
 \be\label{invar}|\Ki(T(z),T(w))|=|\Ki(z,w)|.
 \ee
 A simple way to check \re{invar} is to write $|\Ki(z,w)|$ directly as a function
of the single variable
 $|T_z(w)|$, where $T_z$ is the
isometry taking $z$ to $0$; such an expression is clearly invariant.
It is also possible to get \re{invar} from the invariance of the
process and the covariance formula \re{cov}.

Let $\gamma_{zw}$ be a geodesic connecting $\omega$ and $z$ that
proceeds at speed $s_{zw}=\disti(z,w)$ given by the invariant
distance between $z$ and $w$. Then we have
\begin{eqnarray*}
   | f(z) - f(w) |^2  & = &  \left| \int_0^1 \nablai f(\gamma_{zw}(t) ) \cdot \gamma_{zw}'(t) \, dt \right|^2 \\
   & \le &   s_{zw}^2 \int_0^1 |\nablai f|^2 (\gamma_{zw}(t)) \,
   dt.
\end{eqnarray*}
Thus $\Var_{\rho}(f)$ is bounded above by $|\nablai f|^2$ integrated
against the measure
$$d\vartheta(\zeta)=\frac{1}{2}\int_\Lambda\int_\Lambda \int_0^1
 s_{zw}^2 \delta_{\gamma_{zw}(t)}(\zeta) |\Ki(z,w)|^2 \, dt\,
 \,d\nu(z) \,d\nu(w).
$$
This  is defined in an invariant way, so it must be a some
$\alpha(\rho)\in (0,\infty]$ times the invariant measure. In fact,
`it is the invariant convolution on the symmetric space $\Lambda$
(see, for example \cite{helga} for background) of the standard
invariant measure with the radially symmetric measure
$$
d\vartheta'(\zeta)=\frac{1}{2}\int_\Lambda \int_0^1
 s_{0w}^2 \delta_{\gamma_{0w}(t)}(\zeta) |\Ki(0,w)|^2 \, dt\,
 \,d\nu(w),
$$
and therefore
\begin{eqnarray}\nonumber
\alpha(\rho) =\vartheta'(\Lambda) &=& \frac{1}{2}\int_\Lambda
 s_{0w}^2 |\Ki(0,w)|^2 \,d\nu(w) \int_0^1 \,dt
 \\ &=&
 \pi \int_0^ {R_\Lambda} r s_{0r}^2 |\Ki(0,r)|^2 \eta(r)dr
= \frac{\pi}{\eta(0)} \int_0^{R_\Lambda} r s_{0r}^2 |K(0,r)|^2 dr.
\label{thisq}
\end{eqnarray}
Here $R_\Lambda=\infty$ or $1$ is the radius of the planar model
for $\Lambda$. A direct computation  shows that for all three
models \re{thisq} is bounded by an absolute constant, verifying
our claim.
\end{proof}

\begin{proof}[Proof of Theorem \ref{main}]
Corollary \re{coco} shows that condition CO holds, and the other
conditions have been checked in Section \ref{ouran}. For $f \in C^3$
of compact support Proposition \ref{scum0} gives asymptotic
normality and Lemma \ref{conc} gives the limiting variance, so we
have \be \label{prelim} Z_\rho(f):= \sum_{z\in \mathcal Z} f(z)
 -  \frac{\rho}{\pi}\int_\Lambda |f|d\inv \; \Rightarrow \;{\cal N}
\Bigl( 0,\frac{1}{4 \pi} \|f\|_{H_1(\inv)}^2 \Bigr).
 \ee
Lemma \ref{concbound} allows us to extend the preliminary conclusion
(\ref{prelim}) to the advertised result. For any $f \in
H^1(\Lambda)$ (and of appropriate support) there is a sequence of $
f_{\ep} \in C_0^3$ with $||  f - f_{\ep} ||_{H^1} \ra 0$ as $\ep \ra
0$. Moreover, Lemma \ref{concbound} implies that the family $\{
Z_{\rho} (f) \}$ is tight and also that
\[
   \Bigl|  \ev [e^{i Z_{\rho} (f)} ] - \ev[  e^{i Z_{\rho}  (f_{\ep}) } ] \Bigr|^2 \le \ev \Bigl| Z_{\rho} (f) -  Z_{\rho} (f_{\ep}) \Bigr|^2
                                    \le c \int_{\Lambda}  | \nablai ( f - f_{\ep}) |^2.
\]
The right hand side can be made small at will. Now, choosing a
subsequence $\rho'$ over which $Z_{\rho} (f)$ has a limit in
distribution, we find the Fourier transform of that limit is as
close as we like to that of a mean zero Gaussian with variance
$\frac{1}{4 \pi} \int_{\Lambda} | \nablai f|^2$. (The full limit for
$Z_{\rho} (f_{\ep})$ exists for any $\ep > 0$). Since this appraisal
is the same for any subsequence $\rho'$, we have pinned down the
unique distributional limit of $Z_{\rho}(f)$.
\end{proof}

\bigskip

\noindent {\bf{Acknowledgements.}}  The work of B.R. was supported
in part by NSF grant DMS-0505680,  and that of B.V. by a Sloan
Foundation fellowship, by the Canada Research Chair program, and by
NSERC and Connaught research grants.


\sc \bigskip \noindent Brian Rider, Department of Mathematics,
\\
University of Colorado at Boulder, Boulder, CO 80309.
\\{\tt brider@euclid.colorado.edu},
{\tt www.math.colorado.edu/\~{}brider}

\sc \bigskip \noindent B\'alint Vir\'ag, Departments of Mathematics
and Statistics, \\ University of Toronto, ON, M5S 3G3, Canada.
\\{\tt balint@math.toronto.edu}, {\tt
www.math.toronto.edu/\~{}balint}

\end{document}